\documentclass{amsart}
\usepackage{amsmath}

\begin{document}
\newtheorem{theorem}{Theorem}
\newtheorem{lemma}{Lemma}
\newcommand{\rank}{\mbox{rank}}
\newcommand{\rowspace}{\mbox{rowspace}}
\newcommand{\zfour}{\mathbb{Z}/4\mathbb{Z}}
\newcommand{\kernel}{\mbox{kernel}}
\newcommand{\supp}{\mbox{supp}}
\newcommand{\II}{{\mathcal I}}
\newcommand{\BB}{{\mathcal B}}
\newcommand{\F}{\mathbb{F}}
\newcommand{\C}{\mathbb{C}}
\newcommand{\gram}{\mbox{Gram}}
\newcommand{\ignore}[1]{}
\newcommand{\GF}{\mathrm{GF}}

\title[On the evaluation at $(-\iota,\iota)$ of the Tutte polynomial]{On the evaluation at $(-\iota,\iota)$ of the Tutte polynomial of a binary matroid}

\begin{abstract}
Vertigan has shown that if $M$ is a binary matroid, then $|T_M(-\iota,\iota)|$,  the modulus of the Tutte polynomial of $M$ as evaluated in $(-\iota, \iota)$,  can be expressed in terms of the bicycle dimension of $M$.  
In this paper,  we describe how the argument of the complex number $T_M(-\iota,\iota)$ depends on a certain $\zfour$-valued quadratic form that is canonically associated with $M$. We show how to evaluate $T_M(-\iota,\iota)$ in polynomial time, as well as  the canonical tripartition of $M$ and further related invariants.  
\end{abstract}
\author{R.A.Pendavingh
}
\address{Eindhoven Technical University\\ Den Dolech 2\\ 5600MB Eindhoven\\ the Netherlands}
\email{rudi@win.tue.nl}
\date{\today}
\maketitle
\section{Introduction} 
The {\em Tutte polynomial} of a matroid $M=(E,\II)$ with rank function $r$ is 
\begin{equation}\label{tuttepoly}
T_M(x,y):=\sum_{F\subseteq E} (x-1)^{r(E)-r(F)}(y-1)^{|F|-r(F)}.
\end{equation}
Extending the work of Jaeger, Vertigan and Welsh \cite{JaegerVertiganWelsh1990}, Vertigan investigates the complexity of $\tau^0({\mathcal M}_\F, x,y)$, the problem of evaluating $T_M(x,y)$ given $x,y$ and a matroid $M$ that is linearly represented over $\F$ \cite{Vertigan1998}. Let  $\iota$ denote the imaginary unit, let $\omega:=e^{\iota\frac{2\pi}{3}}$ denote a complex third root, and let $\mathbb{A}$ be the algebraic closure of $\mathbb{Q}$. 
\begin{theorem}[Vertigan, 1998] Let $\F$ be a finite field and $(x,y)\in\mathbb{A}^2$ be a pair other than $(0,0)$, $(1,1)$, and such that $(x-1)(y-1)\neq 1$. Then the problem $\tau^0({\mathcal M}_\F, x,y)$ is \#P-complete, except when 
\begin{enumerate}
\item  $|\F|=2$, and $(x,y)$ is one of $(-1,-1), (0,-1), (-1,0), (\iota,-\iota), (-\iota,\iota)$;
\item $|\F|=3$, and $(x,y)$ is one of $(\omega, \omega^2), (\omega^2, \omega)$; or
\item $|\F|=4$, and $(x,y)$ is $(-1,-1)$.
\end{enumerate} 
\end{theorem}
In the present paper we derive, for a binary matroid $M$, an explicit expression for $T_M(-\iota,\iota)$  that can be evaluated in polynomial time. 

To the best of our knowledge, the complexity status of  $\tau^0({\mathcal M}_\F, 1,1)$ is open. 
Otherwise, the above theorem is now complemented by explicit, polynomial-time computable expressions for the value of the Tutte polynomial on each of the special points:
\begin{enumerate}
\item If $(x,y)=(0,0)$ or  if $(x-1)(y-1)=1$, it is trivial to compute $T_M(x,y)$ for any matroid $M$. 
\item For binary matroids $M$, $|T_M(-\iota,\iota)|$ was determined by Vertigan \cite{Vertigan1998}. Below, we derive an explicit expression for $T_M(-\iota,\iota)$;
\item For ternary matroids $M$, Jaeger \cite{Jaeger1989} has determined $|T_M(\omega,\omega^2)|$, and Gioan and Las Vergnas \cite{GioanLasVergnas2007} found $T_M(\omega,\omega^2)$;
\item For quaternary matroids $M$, $T_M(-1,-1)$ was found by Vertigan \cite{Vertigan1998}, extending a result for graphic/binary matroids by Rosenstiehl and Read \cite{RosenstiehlRead1978}.
\end{enumerate}

The original motivation for this research was a computational problem that arose when writing a matroid package for Sage \cite{sage}. Testing whether two matroids are isomorphic can be made more efficient in practice by comparing matroid invariants, avoiding more involved computation if the values of the invariants do not match. For general matroids, there are few such invariants that are polynomial-time computable. For binary, ternary and quaternary matroids however, the values of  the Tutte polynomial in the above-mentioned special points are clearly the kind of isomorphism invariant we can use for this purpose. With this objective in mind, we shall prove that computing $T_M(-\iota,\iota)$ 
as well as several related invariants takes $O(r(M)^2 |E|)$ time all together.

\section{Preliminaries}
\subsection{Matroids} 
We assume familiarity with matroid theory. In our use of matroid terminology we generally follow Oxley \cite{OxleyBook}. There, a linear matroid $M(A)$ on ground set $E$ is defined from a $k\times E$ matrix $A$. For the present purposes, it will be convenient to also define a linear matroid from a linear subspace, as follows.
 
If $E$ is a finite set, $\F$ is any field and $V$ is a linear subspace of $\F^E$, then $V$ determines a matroid $M(V)$ on $E$ with set of independent sets 
$$\II(V):=\{F\subseteq E\mid \text{there is no }v \in V\text{ such that } \supp(v)\subseteq F\},$$
where $\supp(v):=\{e\in E\mid v_e\neq 0\}$. The relation to the standard definition of linear  matroid is that if $A$ is a $k\times E$ matrix and $V=\ker(A)$, then $M(A)=M(V)$. 

We will denote 
$$V/e:=\{v_{|E-e}\mid v\in V\},$$
so that $M(V/e)=M(V)/e$. 


\subsection{The bicycle dimension}
Let $V$ be a linear subspace of $\F^E$, where $|\F|$ is one of $2,3,4$. The {\em bicycle dimension} of $V$ is defined by
$$d(V):=\dim(V\cap V^\perp).$$
Here $V^\perp:=\{w\in \F^E\mid \langle w,v\rangle =0\text{ for all } v\in V\}$ as usual, where we take $\langle w,v\rangle :=\sum_i w_iv_i$ if $\F=\GF(2)$ or $\GF(3)$  and $\langle w,v\rangle :=\sum_i w_i^*v_i$ if $\F=\GF(4)$, where $^*:\GF(4)\rightarrow \GF(4)$ is the unique nontrivial field automorphism (i.e. $x^*=x^2$ for each $x\in \GF(4)$). 

Vertigan \cite{Vertigan1998} shows that if $V$ is a linear subspace of $\F^E$ and $|\F|$ is one of $2,3,4$, then $d(V)$ depends only on $M(V)$.

\subsection{Quadratic forms}
Let $V$ be a finite-dimensional linear space over  a field $\F$. Let  $b:V\times V\rightarrow \F$ be a bilinear form. 
Then $q:V\rightarrow \F$ is a {\em quadratic form associated with $b$} if
$$ q(\lambda x+\mu y)=\lambda^2q(x)+\mu^2q(y)+\lambda\mu b(x,y)$$
for all $\lambda,\mu\in\F$ and $x,y\in V$. 

A bilinear form  is {\em nondegenerate} if there is no $w\in V$ such that $b(v,w)=0$ for all $v\in V$. 
A quadratic form is {\em nonsingular} if the associated bilinear form is nondegenerate.
A basis $v_1,\ldots,v_k$ of $V$ is {\em orthogonal} with respect to $b$ if 
$$b(v_i, v_j)=0\text{ if } i\neq j.$$
The following is well-known.
\begin{lemma} If $b:V\times V\rightarrow\F$ is a nondegenerate bilinear form and the field $\F$ has characteristic other than 2, then there exists a basis of $V$ that is orthogonal with respect to $b$.\end{lemma}

Two quadratic forms $q,q'$ on $V$ are {\em isomorphic} if there is some linear bijection $L:V\rightarrow V$ such that $q(v)=q'(L(v))$ for all $v\in V$.
If $v_1, \ldots, v_k$ and $w_1,\ldots, w_k$ are $b$-orthogonal bases of $V$ and $q$ is a quadratic form associated with $b$, then $\prod_i q(v_i)$ is a quadratic residue if and only if $\prod_i q(w_i)$ is a quadratic residue. 
Let $\chi(q)$ be $1$ or $-1$ depending on whether the product is a quadratic residue or nonresidue. 
\begin{theorem} If  $q, q'$ are nonsingular quadratic forms on $V$ over a field of characteristic other than 2, then $q$ is isomorphic to $q'$ if and only if $\chi(q)=\chi(q')$.
\end{theorem}
The case when the characteristic is 2 is somewhat more involved. 
A basis $v_1,\ldots,v_{2m}$ of $V$ is {\em alternating} with respect to $b$ if
 $$b(v_i, v_j)=\left\{\begin{array}{ll} 1& \text{if }\{i,j\}=\{k,k+m\}\text{ for some }k\in \{1,\ldots, m\}\\ 0 &\text{otherwise}\end{array}\right.$$
\begin{lemma} If $b:V\times V\rightarrow\F$ is a nondegenerate bilinear form and the field $\F$ has characteristic 2, then exactly one of the following holds:
\begin{enumerate}
\item $V$ has a basis that is orthogonal with respect to $b$, or
\item $V$ has a basis that is alternating with respect to $b$.  
\end{enumerate}
\end{lemma}
Brown \cite{Brown1972} generalizes quadratic forms over $\GF(2)$ to {\em $\zfour$-valued quadratic forms} $q: V\mapsto\zfour$, satisfying
$$ q(x+y)=q(x)+q(y)+\alpha(b(x,y)),$$
where $b:V \times V\mapsto \GF(2)$ is a bilinear mapping, and $\alpha:\GF(2)\rightarrow \zfour$ is the additive group homomorphism such that $\alpha(0)=0$ and $\alpha(1)=2$. 
Such a quadratic form $q$ is {\em nondegenerate} resp. {\em alternating} if and only if the associated bilinear function $b$ is.

Brown also defines an invariant $\sigma(q)$ such that 
\begin{equation} \label{brown}\sum_{x\in V} \iota^{q(x)}= \sqrt{2}^{\dim(V)} e^{\frac{\pi \iota \sigma(q)}{4} }.\end{equation}
Wood \cite{Wood1993} has classified the $\zfour$-valued quadratic forms as follows:
\begin{theorem}[Wood,1993] If $q,q'$ are nonsingular $\zfour$-valued quadratic forms on $V$, then $q$ is isomorphic to $q'$ if and only if  
\begin{enumerate}
\item $\sigma(q)=\sigma(q')$, and 
\item $q$ is alternating $\Leftrightarrow$ $q'$ is alternating.
\end{enumerate}

\end{theorem}

\section{A special point of the Tutte polynomial}
\subsection{Characterization of $T_M(-\iota,\iota)$} Let $V$ be a linear subspace of $\GF(2)^E$, and let $q_V: V\rightarrow \zfour$ be defined by
$$ q_V(x):=| \supp(x)| \mod 4$$
for all $x\in V$. Then $q_V(x+y)=q_V(x)+q_V(y)+\alpha(b(x,y))$ taking $b(x,y)=\sum_i x_iy_i$, so that $q_V$ is a $\zfour$-valued quadratic form on $V$.

If $q_V( y)=0$ for all $y\in V\cap V^\perp$, then for all $x\in V$ and all $y\in V\cap V^\perp$, we have
$$q_V(x+y)=q_V(x)+q_V(y)+\alpha(b(x,y))=q_V(x)$$
and then, we may define $\tilde{q}_V: V/(V\cap V^\perp)\rightarrow \GF(2)$ by setting 
\begin{equation}\label{qtilde}\tilde{q}_V(x + V\cap V^\perp)=q_V(x).\end{equation}
Then $\tilde{q}$ is nonsingular by construction.

We have arrived at the main result of this paper, a characterization of $T_M(-\iota,\iota)$ for binary matroids $M$ in terms of the bicycle dimension and Brown's invariant.
\begin{theorem} Let $V$ be a linear subspace of $\GF(2)^E$ and let $M:=M(V)$. Then
\begin{equation}\label{Tpoint}T_M(-\iota,\iota)= e^{\frac{\pi \iota}{4}(\sigma(\tilde{q}_V)+|E|-3r(E))}\sqrt{2}^{d(V)}\end{equation}
if $q_V(x)=0$ for all $x\in V\cap V^\perp$, and $T_M(-\iota,\iota)= 0$ otherwise.
\end{theorem}
\proof 
By an application of  Greene's formula \cite{Greene1976} (as in \cite[p.390]{Vertigan1998}), we have
$$
\sum_{x\in V} \iota^{q_V(x)}= \iota^{r(E)}(1-\iota)^{|E|-r(E)}T_M(-\iota, \iota).
$$
Rewriting, we obtain
\begin{equation}\label{vertigan}
T_M(-\iota, \iota)= \sqrt{2}^{-|E|+r(E)} e^{\frac{\pi \iota}{4}(|E|-3r(E))}\sum_{x\in V} \iota^{q_V(x)}.
\end{equation}

If  $y\in V\cap V^\perp$, then $V=V+y$ and $b(x,y)=0$ for all $x\in V$, so that $q_V(x+y)=q_V(x)$ for all $x\in V$. Hence if $q_V(y)=2$ for such an $y\in V\cap V^\perp$, then 
$$\sum_{x\in V} \iota^{q_V(x)}=\sum_{x\in V} \iota^{q_V(x+y)}=\iota^{q_V(y)}\sum_{x\in V} \iota^{q_V(x)}=-\sum_{x\in V} \iota^{q_V(x)}.$$
It follows that then $\sum_{x\in V} \iota^{q_V(x)}=0$ and hence $T_M(-\iota,\iota)=0$. 

If on the other hand $q_V(x)=0$ for all $x\in V\cap V^\perp$, then by \eqref{qtilde} we have
$$\sum_{x\in V} \iota^{q_V(x)}=\sum_{w\in V\cap V^\perp} \sum_{v\in \tilde{V}}  \iota^{q_V(v+w)}=\sum_{w\in V\cap V^\perp} \sum_{v\in \tilde{V}}  \iota^{\tilde{q}_V(v+V\cap V^\perp)}$$
where $\tilde{V}$ is any subspace of $V$ so that $V=(V\cap V^\perp) \oplus \tilde{V}$. 
The summation over $V\cap V^\perp$ amounts to a factor $2^{d(V)}$, and using \eqref{brown} on the non-degenerate form $\tilde{q}_V$ we obtain
$$\sum_{x\in V} \iota^{q_V(x)}=2^{d(V)}\sqrt{2}^{\dim(\tilde{V})} e^{\frac{\pi \iota}{4} \sigma(\tilde{q}_V)}.$$
Substituting this expression in \eqref{vertigan}, and using that $$|E|-r(M)=\dim(V)=d(V)+\dim(\tilde{V}),$$ we obtain \eqref{Tpoint}.
\endproof

For comparison, we state the characterization of $T(\omega,\omega^2)$ for ternary matroids due to Gioan and Las Vergnas \cite{GioanLasVergnas2007} in similar terms. For a subspace $V\subseteq \GF(3)^E$, let $q_V:V\rightarrow \GF(3)$ be defined by
$$q_V:x\mapsto |\supp(x)|\mod 3,$$
and let $\tilde{q}_V:V/(V\cap V^\perp)\rightarrow \GF(3)$ be defined by $\tilde{q}_V(x+V\cap V^\perp)=q_V(x)$.
\begin{theorem} Let $V\subseteq \GF(3)^E$ be a linear subspace and $M:=M(V)$. Then 
$$T_M(\omega,\omega^2)=(-1)^{\frac{1-\chi(\tilde{q}_V)}{2}}\omega^{2|E|-r(M)}(\iota\sqrt{3})^{d(V)}.$$\end{theorem}

\subsection{Complexity of computing $T_M(-\iota,\iota)$}
In what follows, let $V\subseteq \GF(2)^E$ be a linear subspace of dimension $k$. A {\em $q$-basis} is a basis $v_1,\ldots, v_k$ of $V$ such that
\begin{enumerate}
\item $v_1,\ldots, v_{k-d(V)}$ is an orthogonal or alternating basis of some subspace $\tilde{V}$ such that $V=\tilde{V}\oplus (V\cap V^\perp)$;
\item $v_{k-d(V)+1},\ldots, v_k$ is a basis of $V\cap V^\perp$; and
\item $q_V(v_{k-d(V)+1})=0,\ldots, q_V(v_{k-1})=0$.
\end{enumerate}
Standard linear algebra techniques yield:
\begin{lemma} \label{qbasis} Given any basis of $V$, computing a $q$-basis  takes $O(\dim(V)^2|E|)$ time. 
\end{lemma}
The following is straightforward from the definition of $q$-basis and \eqref{brown}.
\begin{lemma} Let $v_1,\ldots,v_k$ be a $q$-basis of $V$. Then $q_V(v)=0$ for all $v\in V\cap V^\perp$ if and only if $d(V)=0$ or $q_V(v_{k})=0$. If so, then
\begin{enumerate} 
\item if  $v_1,\ldots, v_{k-d(V)}$ is orthogonal, then 
$$\sigma(\tilde{q}_V)=\#\{i\mid q_V(v_i)=1\} - \#\{i\mid q_V(v_i)=3\} \mod 8,$$
and 
\item if  $v_1,\ldots, v_{k-d(V)}$ is alternating, then 
$$\sigma(\tilde{q}_V) = 4 \#\{i\in \{1,\ldots, m\}\mid q_V(v_{i})=q_V(v_{i+m})\} \mod 8,$$
where $m=\frac{k-d(V)}{2}.$
\end{enumerate}
\end{lemma}

\begin{theorem} Let $V\subseteq \GF(2)^E$ be a linear subspace. Given any basis of $V$, the evaluation of  $T_{M(V)}(-\iota,\iota)$  takes $O(\dim(V)^2|E|)$ time.
\end{theorem}
\proof To compute $T_{M(V)}(-\iota,\iota)$, it suffices to determine $|E|$, $\rank(M(V))=|E|-\dim(V)$, whether $q_V(v)=0$ for all $v\in V\cap V^\perp$, and if so, $\sigma(\tilde{q}_V)$. Given a $q$-basis of $V$, this takes in $O(\dim(V)|E|)$ time.\endproof
As $\dim(V)=|E|-r(M(V))=r^*(M(V))$, this amounts to a complexity bound of  $O(r^*(M(V))^2|E|)$ for evaluating $T_{M(V)}(-\iota,\iota)$ from a basis of $V$.
We note that as in general $T_M(x,y)=T_{M^*}(y,x)=\overline{T_{M^*}(\overline{y}, \overline{x})}$, we have 
$$T_{M(V)}(-\iota,\iota)=\overline{T_{M(V^\perp)}(-\iota,\iota)}$$
We may determine the latter in $O(\dim(V^\perp)^2|E|)=O(r(M(V))^2|E|)$ time from any basis of $V^\perp$.

\subsection{Computing the canonical tripartition}   Let $V\subseteq \GF(2)^E$. 
We consider $$F_i:=\{e\in E\mid d(V/e)=d(V)+i\}.$$
\begin{lemma} \label{canonical} Let $v_1,\ldots,v_k$ be a $q$-basis of $V$. Then
\begin{enumerate}
\item $F_{-1} =\bigcup_{i=k-d(V)+1}^k \supp(v_i)$;
\item $F_1 = \supp(\sum_{i=1}^{k-d(V)}v_i)\setminus F_{-1}$ if the $q$-basis is orthogonal, $F_1=\emptyset$ otherwise; and 
\item $F_0=E\setminus (F_{-1}\cup F_1)$.
\end{enumerate}
In particular, $F_i=\emptyset$ for any $i\not\in \{-1,0,1\}$.
\end{lemma} 
In what follows, we write $A[X,Y]$ for the restriction of a matrix $A$ to the rows indexed by $X$ and the columns indexed by $Y$.
\proof Let  $A$ be any $k\times E$ matrix over $\GF(2)$ such that $V=\rowspace(A)$. Then $V^\perp=\kernel(A)$, and 
$d(V)=\dim(V\cap V^\perp)=k-\rank(AA^T)$. If we write 
$$A^e:=A[\{1,\dots, k\}, E-e]\text{ and }a^e:=A[\{1,\dots, k\}, e]$$ for $e\in E$, then
 $V/e=\rowspace(A^e)$ and $d(V/e)=k-\rank(A^e(A^e)^T)$. Thus
$$d(V/e)-d(V)=\rank(AA^T)-\rank(AA^T+a^e(a^e)^T).$$
Now consider the matrix $A$ whose rows are the given $q$-basis of $V$. 

If $e\in\supp(v_i)$ for some $i>k-d(V)$, so $a^e_i\neq 0$, then consider the matrix $B$ that arises by adding the $i$-th row of $A$ to the rows $j\in \supp(a^e)-\{i\}$. Then the rows of $B$ again form a $q$-basis of $V$, $b^e$ is a unit vector, and hence $d(V/e)-d(V)=\rank(BB^T)-\rank(BB^T+b^e(b^e)^T)=-1$. 

If $a^e_i=0$ for all $i>k-d(V)$, then  $d(V/e)-d(V)=\rank(AA^T)-\rank(AA^T+a^e(a^e)^T)=\rank(BB^T)-\rank(BB^T+b^e(b^e)^T)$ where $B=A[\{1,\ldots, k-d(V)\}, E]$. So without loss of generality, we may assume $d(V)=0$. 

In case the rows of $A$ are orthogonal, let $B:=A[\supp(a^e), E]$. Then $d(V/e)-d(V)=\rank(AA^T)-\rank(AA^T+a^e(a^e)^T)=\rank(BB^T)-\rank(BB^T+b^e(b^e)^T)=\rank(I)-\rank(I+J),$ where $J$ denotes the square all-one matrix with rows and columns indexed by $\supp(a^e)$. 
Then the rank of $I+J$ depends only on the parity of $|\supp(a^e)|$, and we have
$$d(V/e)-d(V)=\left\{\begin{array}{ll}1&\text{if } |\supp(a^e)| \text{ is even}\\
							0&\text{otherwise}\end{array}\right.$$

In case the rows of $A$ are alternating, it remains to show that $d(V/e)-d(V)=0$, i.e. that $AA^T+a^e(a^e)^T$ is nonsingular. If not, there is a nonzero vector $x\in \GF(2)^k$ so that 
$(AA^T+a^e(a^e)^T)x=0$, or equivalently, $AA^Tx=a^e(a^e)^Tx$. As $AA^Tx\neq 0$, we must have $(a^e)^Tx=1$ and hence $AA^Tx=a^e$.  As 
$$AA^Tx=\left[\begin{array}{cc}0&I\\I&0\end{array}\right]x=(x_{k/2+1}, \ldots, x_k, x_1,\ldots, x_{k/2})^T,$$
this implies that $(a^e)^Tx=0$, a contradiction.
\endproof
The triple $(F_{-1}, F_0, F_1)$ is known as the {\em canonical tripartition} of $V$. Originally this notion was developed for graphs by Rosenstiehl and Read \cite{RosenstiehlRead1978}, but the generalization to binary spaces is straightforward. As is clear from the lemma, the determination of the canonical tripartition takes $O(\dim(V)|E|)$ time given a $q$-basis of $V$.

\subsection{A projection} We describe another invariant of subspaces $V\subseteq \GF(2)^E$. The description of this invariant is simpler in the case that $V\cap V^\perp$ is trivial. Therefore, we consider this more restricted setting first.
Following Godsil and Royle \cite{GodsilRoyle2001}, we call a subspace $V\subseteq \GF(2)^E$ {\em pedestrian} if the `bicycle space' $V\cap V^\perp$ is trivial, i.e. if $d(V)=0$. 

Pedestrian spaces are not rare. Fixing a subspace $W\subseteq \GF(2)^E$ of dimension $k$, there are exactly $2^k$ distinct subspaces $V\subseteq \GF(2)^{E+e}$ such that $\dim(V)=k$ and $V/e=W$. Of these spaces $V$, exactly $2^k-2^{k-d(W)}$ satisfy $d(V)=d(W)-1$. From the remaining $2^{k-d(V)}$, exactly half has $d(V)=d(W)$, half has $d(V)=d(W)+1$. A straightforward analysis  shows that if $|E|$ tends to infinity, slightly less than $3/7$-th of the subspaces $V\subseteq \GF(2)^E$ of dimension $k$ will be pedestrian.

If $V\subseteq \GF(2)^E$ is a pedestrian subspace then $\GF(2)^E=V\oplus V^\perp$,  i.e. any vector $x\in \GF(2)^E$ can be uniquely written as $x=v+w$, where $v\in V$ and $w\in V^\perp$. Hence there is a unique linear map $\pi_V: \GF(2)^E\rightarrow V$ so that $\pi_V(x)\in V$ and $x-\pi_V(x)\in V^\perp$ for all $x\in \GF(2)^E$. The matrix $Q_V$ so that $\pi_V(x)=Q_Vx$ is textbook material in linear algebra; it is 
\begin{equation} \label{Q}Q_V=A^T(AA^T)^{-1}A\end{equation}
where $A$ is any matrix with independent rows such that $\rowspace(A)=V$. That $AA^T$ is invertible follows from our assumption that $V$ is pedestrian, that $Q_V\in \rowspace(A)=V$ is clear, and we have $x-Q_Vx\in \kernel(A)=V^\perp$ since
$$ A(x-Q_Vx)=A(x-A^T(AA^T)^{-1}Ax)=Ax-AA^T(AA^T)^{-1}Ax=0.$$
It is not difficult to determine $Q_V$ from a $q$-basis of $V$.
\begin{lemma}Let $V\subseteq \GF(2)^E$ be a pedestrian linear subspace, and let $v_1,\ldots, v_k$ be a $q$-basis of $V$. Then
$$Q_V=\sum_{i=1}^k v_iv_i^T\text{ or } Q_V=\sum_{i=1}^{k/2} v_iv_{i+k/2}^T+v_{i+k/2}v_i^T$$
if the basis is orthogonal or alternating, respectively.
\end{lemma}
\proof Let $A$ be the matrix whose rows are $v_1,\ldots, v_k$. Then $$AA^T=I\text{ or }AA^T=\left[\begin{array}{cc}0&I\\I&0\end{array}\right]$$
depending on whether the $q$-basis is orthogonal or alternating. The result now follows from \eqref{Q}.\endproof
Let $G_V$ be the support graph of the $E\times E$ matrix $Q_V$, so $V(G_V)=E$ and 
$$E(G_V)=\{ef\mid (Q_V)_{ef}\neq 0\}.$$
\begin{theorem} \label{projection}If  $V\subseteq \GF(2)^E$ and $V'\subseteq \GF(2)^{E'}$ are pedestrian binary spaces, then $M(V)$ is isomorphic to $M(V')$ if and only if $G_V$ is isomorphic to $G_{V'}$.\end{theorem}
Binary matroid isomorphism may thus be reduced to graph isomorphism in the time it takes to construct $G_V$, which is $O(\dim(V)|E|^2)$.

If we do not assume that $V$ is pedestrian, then we may still write $$(V+V^\perp)/(V\cap V^\perp)=V/(V \cap V^\perp) \oplus V^\perp/(V\cap V^\perp).$$
Hence there is a unique linear map $\tilde{\pi}_V:(V+V^\perp)/(V\cap V^\perp)\rightarrow V/(V\cap V^\perp)$ such that $$\tilde{\pi}_V(x)\in V/(V\cap V^\perp)\text{ and }x-\tilde{\pi}_V(x)\in V^\perp/(V\cap V^\perp)$$
  for all $x\in (V+V^\perp)/(V\cap V^\perp)$. Clearly, $\tilde{\pi}_V$ coincides with $\pi_V$ if $V$ is pedestrian. 

\begin{lemma} Let $V\subseteq \GF(2)^E$ be a linear space and let $\tilde{V}$ be any subspace of $V$ such that $V=\tilde{V}+(V\cap V^\perp)$. Then $\tilde{\pi}_V(x+(V\cap V^\perp))=\pi_{\tilde{V}}(x)+(V\cap V^\perp)$ for any $x\in V+V^\perp$.
\end{lemma}
\proof  Let $x\in V+V^\perp$. If $\tilde{W}\subseteq V^\perp$ is such that $V^\perp=\tilde{W}\oplus(V\cap V^\perp)$, then $V+V^\perp=\tilde{V}\oplus\tilde{W}\oplus(V\cap V^\perp)$ and we may write $x=v+w+z$, where $v\in \tilde{V}$, $w\in \tilde{W}$, $z\in V\cap V^\perp$. Then 
$$\tilde{\pi}_V(x+(V+V^\perp))=v+(V\cap V^\perp)=\pi_{\tilde{V}}(x)+(V\cap V^\perp)$$
as required.
\endproof
And so we have 
$$ \tilde{\pi}_V(x+(V\cap V^\perp))=Q_{\tilde{V}}x+(V\cap V^\perp)$$
for any $\tilde{V}$ such that $V=\tilde{V}+(V\cap V^\perp)$. Such a $\tilde{V}$ is determined as the span of the first $k-d(V)$ vectors of a $q$-basis  of $V$. 

If $(F_{-1}, F_0, F_1)$ is the canonical tripartition of $V$, then from Lemma \ref{canonical} we have
$$ F_{-1}=\bigcup \{\supp(v)\mid  v\in V\cap V^\perp\}.$$
Then $(Q_{\tilde{V}})_{ef}$ is independent of the choice $\tilde{V}$ if $e,f\in F:=E\setminus F_{-1}$, since then
$$\langle e_e+(V\cap V^\perp), \tilde{\pi}_V(e_f+(V\cap V^\perp))\rangle =e_e^TQ_{\tilde{V}}e_f=(Q_{\tilde{V}})_{ef}.$$
Let $\tilde{G}_V$ denote the support graph of $Q_{\tilde{V}}[F,F]$. Then $\tilde{G}_V=G_V$ if $V$ is pedestrian. Theorem \ref{projection} thus extends to non-pedestrian spaces in a weaker form:
\begin{theorem} If  $V\subseteq \GF(2)^E$ and $V'\subseteq \GF(2)^{E'}$ are binary spaces, then $M(V)$ is isomorphic to $M(V')$ only if $\tilde{G}_V$ is isomorphic to $\tilde{G}_{V'}$.\end{theorem}

\section{Some remarks and conjectures}
\subsection{Counting bases} As was mentioned in the introduction, the complexity of $\tau^0({\mathcal M}_{\F}, 1,1)$  appears to be open. Vertigan announced  that this problem is \#P-complete  for any fixed field $\F$ in \cite{Vertigan1998}, but his result remains unpublished.

It is straightforward from the definition of the Tutte-polynomial \eqref{tuttepoly} that 
$$T_M(1,1)=|\{F\subseteq E\mid r_M(E)=r_M(F)=|F|\}|,$$
i.e. that $T_M(1,1)$ equals the number of bases of $M$. Thus for any fixed field $\F$, the complexity of the following problem is open:
\begin{tabbing}
{\bf given:} \= a basis of a subspace $V\subseteq \F^E$\\
{\bf find:} \> the number of bases of $M(V)$
\end{tabbing}

For regular matroids however,  counting  bases is easy: if $A$ is a totally unimodular matrix with independent rows, then  the number of bases of $M(A)$ equals $\det(AA^T)$. For the case of graphs, this is known as Kirchhoff's Matrix-Tree Theorem --- the proof is a direct application of the Cauchy-Binet Formula for $\det(AB)$ in linear algebra.
Thus computing the number of bases of $M(A)$ clearly takes polynomial time given the totally unimodular matrix $A$. This result generalizes to sixth-root-of-unity matroids and even quaternionic-unimodular matroids. See \cite{PendavinghVanZwam2011} for the extension of the Matrix-Tree Theorem to quaternionic matroids, as well as a noteworthy use of the matrix $Q_V$ in relation to counting bases in minors of $M(V)$. We conjecture that the following related problem may also be solved in polynomial time:
\begin{tabbing}
{\bf given:} \= a dyadic matrix $A$\\
{\bf find:} \> the number of bases of $M(A)$
\end{tabbing}
Here, a rational matrix $A$ is {\em dyadic} if $\det(B)\in \{(-1)^a 2^b\mid a,b\in \mathbb{Z}\}\cup\{0\}$ for each square submatrix $B$ of $A$, and a matroid $M$ is {\em dyadic} if $M=M(A)$ for some dyadic matrix $A$. 

The class of matroids representable over a fixed field has exponential growth rate, but the regular matroids, the sixth-root-of-unity matroids and the dyadic matroids each have a quadratic growth rate (see \cite{GeelenKungWhittle2009}). A bold conjecture one might make is that any minor-closed class of matroids of  {\em quadratic growth rate}  is such that the number of bases of a matroid in the class is computable in polynomial time from some succinct description of the matroid. To be more specific about the nature of this succint description, we conjecture: if $\mathbb{P}$ is a partial field (see \cite{PendavinghVanZwam2011} for a definition) so that the class of matroids representable over $\mathbb{P}$ has quadratic growth rate, then it takes polynomial time to compute the number of bases of $M(A)$, given any $\mathbb{P}$-matrix $A$.

\subsection{Isomorphism testing for binary matroids} We consider the problem:
\begin{tabbing}
{\bf given: } \mbox{ }\= bases for subspaces $V\subseteq \GF(2)^E, V'\subseteq \GF(2)^{E'}$\\
{\bf decide:} \> if $M(V)$ is isomorphic to $M(V')$
\end{tabbing}
This problem properly contains isomorphism testing for 3-connected graphs, and by a simple reduction the general graph isomorphism problem. The complexity of the latter problem remains open to this day. In practice, on may use  a canonical labelling algorithm as was described and implemented by McKay \cite{McKay1981}  for solving such graph isomorphism problems.

We described three isomorphism invariants that one could compute from a basis of $V$ in just $O(\dim(V)|E|^2)$ time: 
\begin{enumerate}
\item $T_{M(V)}(-\iota,\iota)$;
\item the cardinalities of the $F_i$ in the canonical tripartition $(F_{-1}, F_0, F_1)$; and
\item the number of edges of $\tilde{G}_V, \tilde{G}[F_0]$ and $\tilde{G}[F_1]$.
\end{enumerate}
In a random selection of 10,000 subspaces of $\GF(2)^{30}$ of dimension 10, all but 324 of the pairs of were revealed as non-isomorphic by a comparison of these invariants, or were identified as isomorphic by an application of Theorem \ref{projection}. This means that in less than 1 in 100,000 cases, it was necessary to revert to other methods for testing isomorphism. In a forthcoming Sage package for matroid computation, this technique has been implemented to speed up the isomorphism test for binary matroids, and similar methods have been implemented for ternary and quaternary matroids. 

Haggard, Pearce and Royle \cite{HaggardPearceRoyle2010} describe a practical algorithm to compute the Tutte polynomial of a graph, which makes extensive use of graph isomorphism testing to reduce the overall computational effort. A possible application of our isomorphism test would be the extension of their method to an algorithm for computing the Tutte polynomial of binary matroids of moderate size, where the more straightforward methods can only handle small matroids.

\section{Acknowledgements}
We thank Judith Keijsper for a stimulating conversation on the nature of orthogonal projection in binary spaces. We also thank two anonymous referees for their constructive comments.
\bibliographystyle{plain}
\bibliography{math}

\end{document}